\documentclass[preprint,12pt]{elsarticle}
\usepackage{graphicx}
\usepackage{mathtools}
\usepackage{caption}
\usepackage{microtype}
\usepackage[utf8]{inputenc}
\usepackage{hyperref}
\usepackage[english]{babel}
\usepackage{amssymb}
\usepackage{amsfonts}
\usepackage{multirow}
\usepackage{tabularx}

\usepackage{amssymb}
\usepackage{numberbysection}
\begin{document}

\begin{frontmatter}



\title{A Simulation--Based Optimization approach for analyzing the ambulance diversion phenomenon in an Emergency Department network}


\author[lab1]{Christian Piermarini}

\address[lab1]{Dipartimento di Ingegneria Informatica, Automatica e Gestionale  ``A. Ruberti'' \\
              SAPIENZA, Universit\`a di Roma \\ via Ariosto, 25 -- 00185 Roma, Italy \\
              {piermarini.1766526@studenti.uniroma1.it,  roma@diag.uniroma1.it}
              }

\author[lab1]{Massimo Roma}


\begin{abstract}
Most of the studies dealing with the increasing and well--known problem of Emergency Department (ED) overcrowding usually mainly focus on modeling the patient flow within a single ED, without considering the possibilities offered by the cooperation among EDs. Instead, it is important to analyze the overcrowding phenomenon considering an ED network rather than a single ED. 
In this paper, the Simulation--Based Optimization approach is adopted for studying an ED network under different conditions, by means of Discrete Event Simulation models.
In particular we consider the so called \textit{Ambulance Diversion} problem, analyzing different diversion policies.
Such models are carried out from real data collected from six big EDs in the Lazio region of Italy. The aim is to optimize the performances of the entire network, in order to provide the best service to the patients without sustaining too high costs. The obtained experimental results show which are the best diversion policies both in terms of patient waiting time and costs for the service providers.
\end{abstract}

\begin{keyword}
Simulation--based Optimization \sep Emergency Department \sep Discrete Event Simulation \sep Ambulance Diversion.  



\end{keyword}

\end{frontmatter}
\clearpage
\numberbysection
\section{Introduction}
\noindent The Simulation--Based Optimization (SBO) approach is a methodology oriented to the optimization of specific Key Performance Indicators (KPIs) and applicable to a wide range of phenomena and different situations  (see e.g., \cite{gosavi.2014,fu.2015}). It is adopted
to perform deep analyses on management and operative problems
arising in systems which are represented by means of simulation models. The great potentialities of this approach are nowadays well known. Of course, the SBO approach considers Black--Box Optimization problems, since the objective function and the constraints of the optimization problem in hand are not available in analytic form \cite{audet.2017}. As a consequence,
algorithms belonging to the class of Derivative--Free Optimization (DFO) methods (see, e.g., \cite{conn.2009,larson.2019}) must be applied (see also \cite{amaran2016simulation}).
Since both the objective function and the constraints are evaluated from the output of the simulation model, and since the simulation runs could be very time expensive, it is that clear the optimization algorithm adopted should be as effective as possible. Moreover, the interface between the optimization algorithm and the simulation runs should be build so that the communication is handy and not time consuming.
\par
In the framework of emergency medical services, and in particular in tackling ED overcrowding problem, the SBO approach has be successfully applied (see, e.g., the paper \cite{yousefi.2020} where a systematic review of SBO methods applied in hospital EDs is reported and the papers \cite{uriarte.2015,uriarte.2017,chen.2016,feng.2015}). In particular, here we focus on a phenomenon strongly related to the ED overcrowding, namely the \emph{Ambulance Diversion} (AD) (see, e.g. \cite{baek2020centralized,hagtvedt2009cooperative} and the references reported therein). AD is a procedure adopted by ED managers during an emergency situation, when severe patients cannot be accepted and treated anymore. Therefore, incoming ambulances can be diverted to neighboring EDs. More in details, when an ED goes on \emph{AD status}, two plans can be carried out:
 \emph{boarding}, i.e patient joins a queue into the ED hall or inside an ambulance in the parking lot, waiting for his turn; 
\emph{redirection}, i.e. patient is transferred towards another available nearby ED which is not on AD status.
The adoption of the AD represents a great inefficiency of the service provided by an ED. Moreover, this practice is also source of heated discussion: on one hand, the scientific community of physicians does not approve the AD, since every patient should receive the adequate treatment and care whenever it is needed. On the other hand, the ED board and executives aim at obtaining an effective and efficient management of the first--aids, thus possible investments for improving service capacity should always be sustained in the face of an actual request of a specific level of service. As is often the case, the answer lies somewhere in between: even if the AD should not be considered as a standard response to the ED overcrowding, it should not be prohibited when an emergency situation occurs, in particular when the waiting minutes for the interhospital transportation, compared to the hours spent waiting in the ED halls, can make the difference between life or death of a patient.
\par
In this work we formulate the AD problem as an Integer Optimization problem with bound constraints on the decision variables and with simulation constraints. In particular, as case study, we consider a network of six big EDs located in the Lazio region of Italy and we construct a Discrete Event Simulation (DES) model of this network by using ARENA simulation software \cite{ARENA}, The availability of the collected real data enables to build a reliable model. Then an interface written in Python has been created for executing simulation runs via an external code. By using this code we also obtained the communication between ARENA simulation software and an external optimization algorithm and, as far as the authors are aware, this represents a novelty with respect to the standard use of Visual Basic for Application (VBA) included within ARENA software.
This allowed us to use the efficient DFO algorithm for black--box optimization problems with integer variables recently proposed in \cite{liuzzi2020algorithmic}. The obtained results show which is the best AD policy to be adopted in terms of operative performance.
\par
The paper is organized as follows: in Section~\ref{sec:casestudy} the case study is described along with the most commonly used AD policies. Section~\ref{sec:des} reports a description of the DES model we built. In Section~\ref{sec:opt} we describe the formulation of the SBO problem for the case study considered and in Section~\ref{sec:implementation} its implementation. The results of the experimentation we performed is included in Section~\ref{sec:results}. Finally, in Section~\ref{sec:conclusions} some final remarks are reported.

\section{The case study and the different diversion policies}\label{sec:casestudy}
\noindent The purpose of this paper concerns the ED organization and the management, focusing on the well--known problem of the overcrowding. However, this issue is not faced by means of the analysis of a single ED, but considering a network constituted by EDs connected to each other via the transportation on ambulance of severe patients. Thus, the focus is not to manage in detail a single ED in order to optimize some specific KPIs, but to treat an emergency situation widespread in a specific geographic area in order to ensure the safety and the rescue of life--threatening people. 
\par
In particular, the focus is on the \emph{Ambulance Diversion} phenomenon. AD is a procedure put in place by the hospital ED managers when an emergency situation occurs, and severe patients cannot be accepted and treated anymore, usually indicated as \textit{ED on AD status}. Two possibilities are usually considered, with possible different policies adopted, when an ED goes on an AD status:
\begin{itemize}
	\item the \emph{boarding}, namely the patients are anyhow accepted but, due to unavailability of ED resources, they wait (in a possible long queue) for visits and treatments. They are usually placed into and ED hall or inside the ambulance in the parking lot during the waiting time;
	\item the \emph{redirection}, namely the patients are transferred to another available nearby ED.
\end{itemize}
\par
In particular, as case study we consider a network composed by six EDs located in the Lazio region of Italy. For these EDs, data for the whole year 2012 have been collected, only considering yellow and red tagged patients\footnote{We refer to the standard four--colors tag scale used in the ED triage}, since AD only concerns severely damaged or life--threatening patients. 
For each single ED and for each patient the following timestamps are available: the starting time of the triage, the starting time of the first medical visit and the discharge time. Moreover, the times needed for the transportation of patients between each ED of the network is known.
\par 
Another important concept related to the AD phenomenon is the definition of \emph{diversion policies}, which refers to the set of organizational choices that defines unequivocally which are the conditions that establish when an ED must be considered on AD status, if an ED has to execute the patient's boarding or the redirection, and which are the destinations of the redirected patients. Of course, the diversion policies strongly affect the dynamics of the whole ED network. In this paper, we consider the following four most commonly adopted policies:
\begin{enumerate}
	\item[\textbf{P1}:] \textit{Ambulance stoppage}, which does not allow any redirection of patients inside the network, authorizing the boarding instead.
\item[\textbf{P2}:] \textit{Redirect towards the nearest ED (complete diversion)} when all the ED resources are occupied and patients are redirected to the nearest ED, if the latter is not on diversion.
\item[\textbf{P3}:] \textit{Redirect towards the nearest ED, with priority to the red--tagged patients (partial diversion)}, which is identical to the previous one, but it involves the diversion for only yellow--tagged patients. More specifically,  when the number of ED occupied resources reaches a certain threshold value, the yellow--tagged patients are redirected to the nearest ED, if the latter is not on partial or complete diversion.
\item[\textbf{P4}:] \textit{Redirection to the least occupied ED of the network}, regardless of the distance between the ``sending'' and the ``receiving'' ED.
\end{enumerate}

\section{The Discrete Event Simulation model}\label{sec:des}
\begin{figure}[h]
	\centering
	\includegraphics[width=10.5cm]{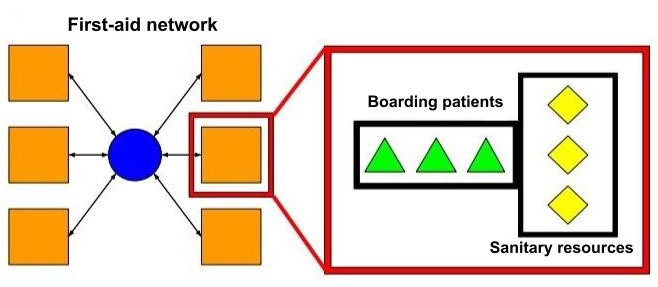}
	\caption{ED network and single ED structure}\label{fig:structure}
\end{figure}

\noindent In order to build the DES model to be used for studying the AD policies for the ED network considered, we first introduce the concept of sanitary resource. We define \textit{sanitary resource} an integer valued parameter which corresponds to the set of all human (like physicians, nurses) and physical (like ED rooms, beds, stretchers) resources required to treat a single patient. The introduction of this notion enables to model each single ED as single queuing system, in which the sanitary resource represents the number of servers as represented in Figure~\ref{fig:structure}
As we will describe afterward, the value of the sanitary resource associated to each ED will be determined by means of a model calibration phase, so that 
this simplified representation of each single ED allows us to obtain a good accuracy in evaluating the KPIs of interest computed from the simulation output.
Note these sanitary resources actually will result as control variables during the optimization phase.
\par
In order to develop the DES model based on the structure now described, it is required to know: 
\begin{itemize}
\item the patient arrival stochastic process to each ED of the network; \item the probability distribution of the Length of Stay (LOS) at each ED of the network. 
\end{itemize}
To this aim, we consider 8-hours time slots of each day (from 00:00 to 08:00, from 08:00 to 16:00 and from 16:00 to 24:00). The choice of 8-hours time slots is motivated by the fact that the ED staff is scheduled on 8-hours shifts, hence an analysis based on these time slots results quite realistic. 
\par
We evaluate patient interarrival time and LOS probability distributions for each time slot and for each yellow and red--tagged patient. The data were collected into Excel files where each line is associated to a patient and reports:
the ED name, the triage tag assigned, date and time of the starting of the triage, date and time of the starting of the visit and date and time of discharge from the ED. Of course, from these data we can easily obtain the time difference between the starting of the triage and the starting of the visit. Since the starting of the visit is the real tacking charge of the patient by the ED, this time difference actually represents the patient waiting time before the ED tacking charge and we will use this as the main KPI to calibrate and validate our DES model.
\par
As regards, the stochastic process of patient arrival at each ED, according to the standard assumption in the literature (see, e.g.,\cite{Whitt.2017,Kuo.2016,Zeinali.2015,Ahmed.2009,Ahalt.2018,Guo.2017,DGLMR:FSMJ}),
we model patient arrivals by a Nonhomogeneous Poisson Process (NHPP) and following a usual procedure we approximate the arrival rate by a piecewise constant function. On the basis of the data, we determine the patient arrival rate for each ED of the network, for each time slot and for each color tag. In Table~\ref{tab:arrivals}, we report the number of patient arrivals during the whole year at each ED of the network, for red and yellow--tagged patients, in each time slot.
\begin{table}[htbp]
	\centering
\begin{tabular}{r|c|ccc|c|}
	& \textit{Color tag} & \footnotesize{00:00-08:00} & \footnotesize{08:00-16:00} & \footnotesize{16:00-24:00} & \textit{Total} \\ \hline
\multirow{2}*{ED1} &  Yellow & 875 & 2688 & 2279 & 5842 \\
& Red & 94 & 299 & 187 & 580 \\ \hline
\multirow{2}*{ED2} 	
&  Yellow & 390 & 1592 & 1299 & 3281 \\
& Red  & 56 & 153 & 107 & 316 \\ \hline
\multirow{2}*{ED3} 	
&  Yellow & 599 & 2331 & 1839 & 4769 \\
& Red & 84 & 102 & 92 & 278 \\ \hline
\multirow{2}*{ED4} 
&  Yellow & 439 & 2118 & 1503 & 4060 \\
& Red & 70 & 295 & 202 & 567 \\ \hline
\multirow{2}*{ED5} 
&  Yellow & 399 & 2007 & 1597 & 4003 \\
& Red & 55 & 299 & 188 & 542 \\ \hline
\multirow{2}*{ED6} 	
&  Yellow & 324 & 1138 & 832 & 2294 \\
& Red & 37 & 58 & 72 & 167 \\ \hline
\end{tabular}
	\caption{Number of patient arrivals at each ED, for each color tag and for each time slot}
	\label{tab:arrivals}
\end{table}
Note that the red--tagged patient arrivals in each time slot (and in particular during the night time slot) are few. This is another motivation for considering 8-hours time slots, since smaller time slots would imply an even smaller number of patient arrivals, making meaningless every statistical analysis.
\par
Regarding the LOS probability distributions, the time interval between the take charge and the relative discharge has been evaluated from available data for each ED of the network, for each patient and for each color tag. Then, a standard statistical analysis has been performed to obtain the best fitting. 
\par
\begin{figure}[htbp]
	\centering
	\includegraphics[height=4truecm,width=\linewidth]{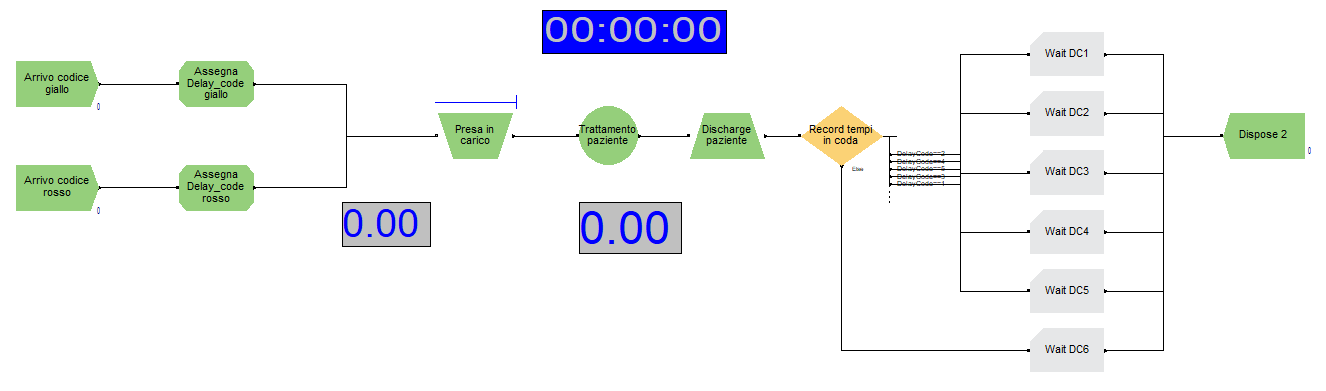}
	\caption{ARENA model of a single ED }\label{fig:singleED}
\end{figure}

The model has been implemented by using {\sf ARENA 16.1} simulation software \cite{ARENA}. On the basis of the assumptions previously described, each ED can be represented as reported in Figure~\ref{fig:singleED}.
The model entities are the patients which flows through the ED; the first two modules on the left create the flow of red and yellow--tagged patients inside the model. An attribute (\emph{delay\_code}) is assigned to each entity, depending on the color tag and the time slot in which the patients entered the ED. Then, the entities join a queue, with priority to the red tagged--patients, and then they are ``taken in charge'' by one free sanitary resource (if any) which is seized to this aim. The patients stay inside the ED for a time that depends on the LOS probability distribution of the associated \emph{delay\_code}. At the end, the entity releases the seized resource and exit the ED. If one sanitary resource is not available, entity is queued waiting his turn. Record modules are used to store the waiting time value of each patient, distinguishing each color tag and each time slot. This waiting time is computed as time difference between the starting of the triage and the actual take charge of the patient (the starting of the visit), and thus it represents the waiting time of each entity before seizing a sanitary resource.
\par
Since the AD phenomenon strictly depends on patient waiting time, we use this KPI to calibrate our model. Namely, we performed an accurate model calibration phase in order to determine the sanitary resource value of each ED so that the average patient waiting time computed from simulation output is a good approximation of the real average patient waiting time determined from the available data. This has been obtained by determining the sanitary resource values which minimize 
$$ \sum_{i=1}^{3} \sum_{j=1}^{2}\ \lvert  W^{sim}_{i,j} - W^{real}_{i,j} \rvert,$$
where $i=1,2,3$ represents the time slot, $j=1,2$ the color tag (red and yellow ones) and $W^{sim}_{i,j}$ and $W^{real}_{i,j}$ are the average patient waiting time obtained from simulation output and from data, respectively. The results of the calibration procedure are reported in Table~\ref{tab:resultscal}.
\begin{table}[htbp]
	\centering
\begin{tabular}{c|ccc|}
		& \footnotesize{00:00--08:00} & \footnotesize{08:00--16:00} & \footnotesize{16:00--24:00} \\ \hline
ED1     & 4 & 4 & 4 \\
ED2     & 4 & 5 & 4 \\
ED3     & 4 & 4 & 3 \\
ED4     & 4 & 5 & 2 \\
ED5     & 4 & 5 & 3 \\
ED6     & 3 & 2 & 2 \\ \hline
\end{tabular}
\caption{Sanitary resource values obtained from the calibration procedure}
\label{tab:resultscal}
\end{table}
By using this values of the sanitary resources, we obtain a model for each ED reproducing its ``as--is'' status, without considering any AD.
\par
\begin{figure}[h]
			\centering
			\includegraphics[width=12cm]{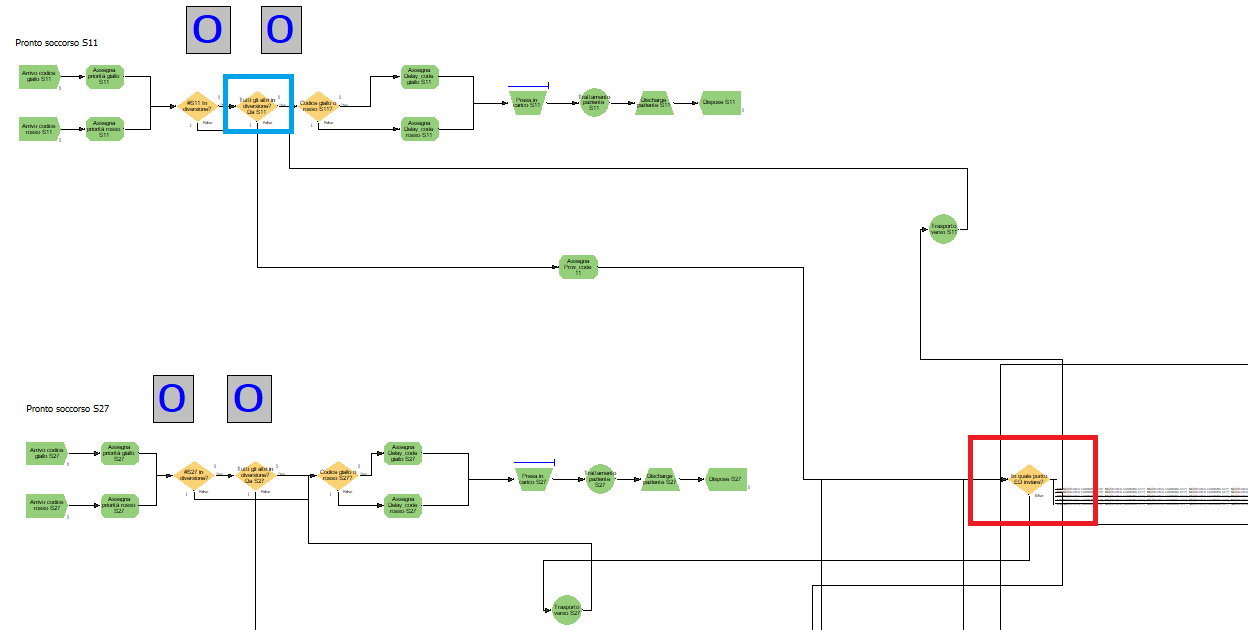}
			\caption{ARENA module (part) for the AD policy {\bf P4}. The red square highlights the central decision module.}\label{fig:P4}
\end{figure}

Then we implemented the ED network composed by the six EDs under study, reproducing the four AD policies (\textbf{P1}, \textbf{P2}, \textbf{P3}, \textbf{P4}) described in Section~\ref{sec:casestudy}. The first one (\textbf{P1}) does not consider any redirection, so that each ED is considered as a single first--aid point.
In the second one (\textbf{P2}), to implement complete redirection to the nearest ED, some decision modules are added to each ED model in order to represent the diversion condition and the actual redirection of the patients when all the sanitary resources are seized and a red/yellow--tagged patient arrives. Policy \textbf{P3} considers partial diversion and hence only redirection of yellow tagged--patients has been implemented, depending on a threshold value of the sanitary resources. Finally, policy \textbf{P4} involves an additional central decision module for redirection towards the ED of the network with the least number of sanitary resources seized. Figure~\ref{fig:P4} report a part of such ARENA model, where the mentioned decision module is highlighted.

\section{Statement of the optimization problem}\label{sec:opt}
\noindent In this section we report the optimization problem we formulate in order to determine which is the AD policy to be considered the most effective to provide prompt first--aid to red and yellow--tagged patients.
The goal is to determine the value of the sanitary resource of each ED and the AD policy to be adopted, so that some suited KPIs are optimized.
\par
We introduce the integer valued decision variables $n_{ij}$ to indicate the value of the sanitary resource assigned to the $i$-th ED ($i=1, \ldots ,6$) during the $j$-th time slot ($j=1, \ldots ,3$). These 18 variables represent the \textit{control variables} of the SBO problem.
\par
As regards the KPIs of interest, we consider the patient average  Non-Valued Added (NVA) time associated to each ED, for red and yellow--tagged patients, corresponding to waiting and transfer times. We denoted by $t^{\rm NVA/R}_i$ and $t^{\rm NVA/Y}_i$ the average NVA time of the $i$-th ED ($i=1, \ldots ,6$) for red and yellow--tagged patients, respectively. It is important to note that such times actually are function of the variables $n_{ij}$ and they can be only computed from the output of the simulation runs. Therefore the following expressions
$$
\sum_{i=1}^6 t^{\rm NVA/R}_i \qquad \hbox{and} \qquad\sum_{i=1}^6 t^{\rm NVA/Y}_i
$$
constitute two important KPIs to be minimized. 
\par
Of course, to keep costs low, the overall number of sanitary resources should be minimized, too. Of course, this is a conflicting objective with respect to the minimization of the NVA times defined above. Therefore, actually, we have a black box multiobjective problem which we tackle by reducing into a single objective problem by means of the standard scalarization procedure (see, e.g., \cite{miettinen:99}). More in detail, we consider the following objective function
\begin{equation}\label{eq:fo}
w_1 \cdot 480 \sum_{i=1}^6 \sum_{j=1}^3 n_{ij}+w_2 \sum_{i=1}^6 t^{\rm NVA/Y}_i + w_3 \sum_{i=1}^6 t^{\rm NVA/R}_i
\end{equation}
where $w_1$, $w_2$ and $w_3$ are suited weights. The last two terms represent the overall NVA times (in minutes) for red an yellow--tagged patients, respectively; the first term measures the overall time (in minutes) for which the sanitary resources are scheduled in the 8 hours (480 minutes) time slots. This latter term takes into account the number of sanitary resources assigned to each ED and it has been expressed in terms of time (in minutes) for homogeneity with respect the other two terms of the objective function.  This term is meant to depict the total cost to be sustained in order to maintain a certain number of sanitary resources active in each time slot of each ED.
%
%After checking the correct behaviour of the models, the optimization phase began. First of all, the formulation of the problem has been carried out. The problem has only one objective function, with eighteen different control variables. The structure of the problem is as follows:
\par
As regards the constraints, we have bound constraints on the variables $n_{ij}$ and simulation constraints which prevent the NVA times to exceed specific threshold values for each color tag. Such values are 40 minutes for the yellow--tagged patients and 20 minutes for the red--tegged ones.
\par
The complete formulation of the SBO problem we consider is the following:
\begin{equation}\label{eq:sboproblem}
\begin{aligned}
	\min \quad & w_1 \cdot 480 \sum_{i=1}^6 \sum_{j=1}^3 n_{ij}+w_2 \sum_{i=1}^6 t^{\rm NVA/Y}_i + w_3 \sum_{i=1}^6 t^{\rm NVA/R}_i \\
	{\rm s.t.} \quad & t^{\rm NVA/Y}_i \leq 40, \qquad i=1, \ldots ,6, \\
	     & t^{\rm NVA/R}_i \leq 20, \qquad i=1, \ldots ,6, \\
	     & 2 \leq n_{ij}\leq 10, \quad ~~ i=1, \ldots ,6; ~j = 1, \ldots ,3, \\ 
	     & n_{ij} \quad \hbox{integer},
\end{aligned}
\end{equation}
which is a Black Box Integer Optimization (BBIO) problem with 18 variables (with the corresponding bound constraints) and 12 simulation constraints.

\section{The SBO implementation}\label{sec:implementation}
\noindent To solve the BBIO Problem~\eqref{eq:sboproblem}, we use an algorithm belonging to the class of Derivative Free Optimization methods recently proposed in \cite{liuzzi2020algorithmic}, named {\sf DFLINT}\footnote{Publicly available in the {\em DFL Library} at the URL {\tt www.iasi.cnr.it/$\sim$liuzzi/DFL}}. It is an algorithm for black box inequality and box constrained integer nonlinear programming problems, hence well suited for solving Problem~\eqref{eq:sboproblem}. Since the {\sf DFLINT} code is available in Python, we decided to code Problem~\eqref{eq:sboproblem} in Python as well. Moreover, we also created a specific interface between the optimization algorithm and the ARENA simulation package by using Python.
Note that ARENA is a Windows desktop application and no information is provided for executing simulation runs from an external software or for
building an interface with any external sources. Only VBA modules can be used within ARENA package and we believe that their use for connecting an external optimization engine with an ARENA model is really cumbersome. As far as the authors are aware, the use of an interface coded in Python for enabling execution of simulation runs and communication between ARENA and external software is a novelty which allowed us to significantly improve the overall performance of the SBO approach applied to the problem under study.
For this reason, we now report a brief description of our implementation, with some technical details. A more insight description is available in \cite{piermarini2021}.
\par
Starting from the values reported in Table~\ref{tab:resultscal} (assumed as the starting point of the optimization algorithm), at each iteration of the algorithm, the values of the control variables (the sanitary resources) are updated and written in an Excel file so that they can be used as input parameters of the simulation. Then simulation run is executed via the Python code. This is performed by using four executable ARENA files {\tt expmt.exe}, {\tt model.exe}, {\tt linker.exe} and {\tt siman.exe} which must be executed in this specific order. After creating the two files of the model (with {\tt .exp} and {\tt .mod} extension), by using the {\tt subprocess.run} Python module, the aforementioned files are converted into two other files (with {\tt .e} and {\tt .m} extensions). Those last two files (in this order) are required for the {\tt linker.exe} to work. As result, a {\tt .p} file is produced and the {\tt siman.exe} program uses it to run a simulation. Finally, the last part of the code,
aims at opening the ARENA output file and extracting the requested values from it.

In the ARENA output file, the NVA times of interest, are stored for every single replication and, based on the Sample Average Approximation approach \cite{pagnoncelli.2009,kim.2015-b}, we compute the average values of the NVA times as mean over independent replications. The corresponding confidence interval is computed, too.

\section{Experimental results}\label{sec:results}
\noindent In this section we report the results obtained by means of the SBO approach we propose, applied to the ED network under study. The aim is to determine which AD policy (among those described in Section~\ref{sec:casestudy}) and which ED settings (in terms of sanitary resources to be allocated) perform the best in the sense that the objective function \eqref{eq:fo} is minimized. In our experimentation the weights $w_1$, $w_2$ and $w_3$ in \eqref{eq:fo} are chosen so that the order of magnitude of the three terms of the objective function \eqref{eq:fo} is the same, namely $w_1=1$, $w_2=300$ and $w_3=600$. Of course different values could be experimented, leading to different weighing of the three terms of the objective function. This shows the flexibility of our approach, since it enables to assign higher weights to terms which are to be considered most relevant, for instance to the term which refers to the most critical (red--tagged) patients.
\par
As regards the simulation runs, we perform 30 independent replications for each run, since we experimented that this choice guarantees a good accuracy of the output results of the simulation. The length of each replication is 365 days and we adopt a warm--up period of 48 hours since this ensures the system to reach the steady state.
\par
As concerns the parameters of the {\sf DFLINT} algorithm, we use default parameters of the original code. Moreover, in order to limit the overall computational burden, the maximum number of function evaluations (which is used as stopping criterion of the algorithm) is set to 700.
\par
In Table~\ref{tab:resultsopt} we report the optimal value of the sanitary resources for each considered policy, for each ED and for each time slot determined by the optimization algorithm. The values corresponding to the starting point (see Table~\ref{tab:resultscal}) are also reported for an easy comparison.
\begin{table}[htbp]
\centering
\begin{tabular}{r|c|c|cccc|}
	     &                    & \textit{Starting} & & & & \\
		 & \textit{Time slot} & \textit{point} & \textbf{P1} & \textbf{P2} & \textbf{P3} & \textbf{P4} \\ \hline
		\multirow{3}*{ED1} &  \footnotesize{00:00--08:00} & 4 & 7 & 4 & 6 & 5 \\
		& \footnotesize{08:00--16:00} & 4 & 5 & 4 & 3 & 5 \\
		& \footnotesize{16:00--24:00} & 4 & 4 & 3 & 5 & 5 \\ \hline
		\multirow{3}*{ED2} &  \footnotesize{00:00--08:00} & 4 & 6 & 5 & 5 & 5 \\
		& \footnotesize{08:00--16:00} & 5 & 5 & 5 & 4 & 5 \\
		& \footnotesize{16:00--24:00} & 4 & 7 & 4 & 4 & 4 \\ \hline
		\multirow{3}*{ED3} &  \footnotesize{00:00--08:00} & 4 & 4 & 3 & 4 & 4 \\
		& \footnotesize{08:00--16:00} & 4 & 6 & 4 & 5 & 4 \\
		& \footnotesize{16:00--24:00} & 4 & 3 & 4 & 3 & 4 \\ \hline
		\multirow{3}*{ED4} &  \footnotesize{00:00--08:00} & 4 & 7 & 5 & 6 & 4 \\
		& \footnotesize{08:00--16:00} & 5 & 4 & 6 & 6 & 3 \\
		& \footnotesize{16:00--24:00} & 2 & 6 & 4 & 4 & 4 \\ \hline
		\multirow{3}*{ED5} &  \footnotesize{00:00--08:00} & 4 & 7 & 7 & 6 & 5 \\
		& \footnotesize{08:00--16:00} & 5 & 4 & 6 & 4 & 3 \\
		& \footnotesize{16:00--24:00} & 2 & 6 & 4 & 5 & 5 \\ \hline
		\multirow{3}*{ED6} &  \footnotesize{00:00--08:00} & 3 & 4 & 4 & 4 & 3 \\
		& \footnotesize{08:00--16:00} & 3 & 4 & 4 & 3 & 3 \\
		& \footnotesize{16:00--24:00} & 2 & 4 & 3 & 3 & 2 \\ \hline
\end{tabular}		
	\caption{Values of the sanitary resources corresponding to the starting point and to the optimal point for each AD policy}
	\label{tab:resultsopt}
\end{table}
It can be easily observed that an increase of the number of sanitary resources (with respect to the starting point) in many cases is required to ensure the viability of the corresponding AD policy. This is mainly due to the satisfaction of the simulation constraints in Problem~\ref{eq:sboproblem} which impose that the NVA times must be below a specific value.
\par
In Table~\ref{tab:resultsfopt} we report the corresponding optimal function values for each AD policy.
\begin{table}[htbp]
	\centering
	\begin{tabular}{r|c|c|}
		& \textit{Objective function value} & \textit{Optimal objective} \\
		& \textit{at the starting point} & \textit{function value} \\ \hline
\textbf{P1}    & 127454,63 & 51975,77  \\
\textbf{P2}     & 92956,13 & 43506,65  \\
\textbf{P3}   & 93758,85 & 46180,20  \\
\textbf{P4}      & 47926,04 & 42037,83  \\ \hline
\end{tabular}		
\caption{Optimal values of the objective function for each AD policy}
\label{tab:resultsfopt}
\end{table}
The obtained improvement is clearly evidenced for all the policies considered, and it is significant in most cases.
\par
To give more complete detailed results, in the Appendix we report the values of all the NVA times corresponding to the starting point (Table~\ref{tab:app1}) and to the optimal solution (Table~\ref{tab:app2}), for each AD policy, for each ED and for each color tag.
\par
By observing the results in Table~\ref{tab:resultsfopt}, the best diversion policy in terms of the objective function \eqref{eq:fo} can be easily determined: it is the policy {\bf P4}, followed by the policy {\bf P2}. The former one, in particular, provides good operative performances even before the optimization, i.e. in correspondence of the starting point. Furthermore, the worst diversion policy is {\bf P1} for which the best improvement of the objective function is registered.
Giving higher priority to the red code while imposing the redirect towards the nearest ED (policy {\bf P3}) seems to worsen the performances of the ED network. The reasons behind this effect can be explained by the analysis of the average NVA times registered for this diversion policy (see
Tables~\ref{tab:app1} and \ref{tab:app2} in the Appendix).
In fact, the average times referred to the red--tagged patients are significantly low when compared to the ones obtained with the other diversion policies, but this leads to an increase of the average NVA times referred to the yellow--tagged patients. On the overall, the result is an overbalance towards a general worsening of the operative performances of the network, which presumably derives from an excessive priority given to the red--tagged patients that have a lower probability of showing up in an ED.
\par
By observing in detail Tables~\ref{tab:app2} in the Appendix corresponding to the obtained optimal point, it is clear that in most cases the NVA times are similar across the different diversion policies.
This is also due to the queue discipline adopted for boarding patients that gives highest priority to red--tagged patients. It can be also noticed that, the average NVA times do not exceed 5 minutes, which can be considered a really good operative result.
\par
Finally, it is possible to understand which is the main issue associated with the ambulance stoppage diversion policy {\bf P1}. Indeed, the inability to redirect patients toward other EDs has the direct consequence of an allocation of an higher number of sanitary resources to maintain good operative performances. However, similar performances, in terms of NVA times, can be obtained by applying other diversion policies with a reduced
number of sanitary resources, thanks to the interhospital transportation of the patients. In other words, without adopting redirection policies, the same average NVA times are obtained but with higher costs due to an increase of the sanitary resources needed. The key element of this consideration is that a good trade--off is needed between the reduction of the patient average waiting time inside the ED and the average transportation time towards another ED. This analysis can be performed by different choices of the weights $w_1$, $w_2$ and $w_3$ in the objection 
function of the SBO Problem \eqref{eq:sboproblem} we formulated. In fact,
this enables to obtain a different balance of the three terms (referred to three single KPIs) that compose the objective function. For instance, an higher cost containment can be obtained by increasing $w_1$, or the same priority level (as regards the AD policies) can be assigned to red and yellow--tagged patients by selecting $w_2=w_3$.

\section{Conclusions}\label{sec:conclusions}
\noindent In this paper we deal with the AD phenomenon by the SBO approach. In particular, we formulate the problem as a Black--Box Integer Optimization problem which we solve by means of an effective DFO algorithm recently proposed in literature. In order to experiment the approach proposed, we consider a network composed by six EDs located in the Lazio region of Italy. The data collected for these EDs (related to one year) have been used for constructing a DES model which we implemented by using ARENA simulation software. A specific interface, written in Python, has been built for connecting the ARENA model and the optimization algorithm. 
\par
The obtained results on the case study considered show that the best diversion policy in terms of operative performance seems to be the redirection towards the least occupied ED. Anyhow, our experimentation clearly evidences that the ambulance diversion and the interhospital transportation of the patients are procedures that should be always considered, especially during important emergency situations. 
\par
As future work, we believe that it would be very interesting to further deepen the analysis on the ambulance diversion policies through dedicated research activities, possibly supported by higher quality information and more efficient data management tools. 
%
%It is also important to remember that the lines of code that have been used for the communication between simulator and optimizer are specifically referred to the ARENA Simulator, so it would be appealing to exploit the structure of the code presented to develop other algorithms compatible with different simulators. 

\appendix
\section{Detailed results on the NVA times}
\noindent
In this appendix we include detailed results concerning the NVA times $t^{\rm NVA/Y}_i$ and $t^{\rm NVA/R}_i$, $i=1, \ldots ,6$ for each AD policies, for each ED and for each color tag. The tables report the average values and the half--width of the $95\%$ confidence interval.
In particular, in Table~\ref{tab:app1} we report the values of $t^{\rm NVA/Y}_i$ and $t^{\rm NVA/R}_i$ corresponding to the \textit{starting point}. In Table~\ref{tab:app2} we report the values of $t^{\rm NVA/Y}_i$ and $t^{\rm NVA/R}_i$ corresponding to the \textit{optimal point}. All the times are expressed in minutes.
\begin{table}[htbp]
	\centering
\begin{tabular}{r|c|ccccc|}
	 & \textit{Color tag} & {\bf P1} & {\bf P2} & {\bf P3} & {\bf P4} \\ \hline
	\multirow{2}*{ED1} 	
	&  Yellow & $15.89\pm 0.96$ & $8.55\pm 0.45$ & $12.78\pm 0.55$ & $2.98\pm 0.09$ \\
	& Red & $6.78\pm 0.44$ & $4.86\pm 0.42$ & $4.19\pm 0.36$ & $2.70\pm 0.18$ \\ 
	\hline
	\multirow{2}*{ED2} 	
	&  Yellow & $25.77\pm 4.15$ & $2.05\pm 0.39$ & $6.89\pm 0.48$ & $2.65\pm 0.17$ \\
	& Red & $11.18\pm 1.62$ & $1.72\pm 0.22$ & $0.74\pm 0.12$ & $2.43\pm 0.30$ \\ 
    \hline
	\multirow{2}*{ED3} 	
	&  Yellow & $5.96\pm 0.97$ & $2.10\pm 0.03$ & $5.06\pm 0.87$ & $1.32\pm 0.17$ \\
	& Red & $3.94\pm 0.74$ & $1.53\pm 0.33$ & $1.98\pm 0.42$ & $1.21\pm 0.28$ \\ 
    \hline
	\multirow{2}*{ED4} 	
	&  Yellow & $50.28\pm 2.75$ & $34.13\pm 1.65$ & $43.43\pm 1.54$ & $5.10\pm 0.27$ \\
	& Red & $24.40\pm 1.08$ & $20.19\pm 1.08$ & $15.83\pm 0.75$ & $4.46\pm 0.34$ \\ 
    \hline
	\multirow{2}*{ED5} 	
	&  Yellow & $54.63\pm 2.98$ & $37.20\pm 1.57$ & $46.48\pm 1.60$ & $5.27\pm 0.32$ \\
	& Red & $23.70\pm 1.43$ & $17.85\pm 1.01$ & $14.22\pm 0.81$ & $4.37\pm 0.31$ \\ 
    \hline
	\multirow{2}*{ED6} 	
	&  Yellow & $12.51\pm 1.26$ & $12.98\pm 1.25$ & $9.05\pm 0.93$ & $2.33\pm 0.19$ \\
	& Red & $7.10\pm 1.22$ & $7.47\pm 0.96$ & $4.66\pm 0.70$ & $2.08\pm 0.41$  \\ 
    \hline
\end{tabular}		
\caption{Values of the NVA times $t^{\rm NVA/Y}_i$ and $t^{\rm NVA/R}_i$, $i=1, \ldots ,6$ corresponding to the starting point}
\label{tab:app1}
\end{table}		
\newpage
\begin{table}[htbp]
	\centering
\begin{tabular}{r|c|ccccc|}
	& \textit{Color tag} & {\bf P1} & {\bf P2} & {\bf P3} & {\bf P4} \\ \hline
\multirow{2}*{ED1} 
&  Yellow & $1.29\pm 0.29$ & $2.39\pm 0.14$ & $3.23\pm 0.15$ & $1.18\pm 0.05$ \\
& Red & $0.87\pm 0.14$ & $1.97\pm 0.16$ & $0.40\pm 0.12$ & $1.16\pm 0.10$ \\ 
\hline
\multirow{2}*{ED2} 
&  Yellow & $2.62\pm 0.43$ & $1.25\pm 0.17$ & $4.88\pm 0.16$ & $1.18\pm 0.11$ \\
& Red & $1.54\pm 0.33$ & $1.07\pm 0.22$ & $0.33\pm 0.05$ & $1.09\pm 0.16$ \\ 
\hline
\multirow{2}*{ED3} 
&  Yellow & $0.99\pm 0.23$ & $1.22\pm 0.20$ & $1.78\pm 0.07$ & $0.84\pm 0.07$ \\
& Red & $0.65\pm 0.20$ & $1.08\pm 0.23$ & $0.11\pm 0.05$ & $0.80\pm 0.11$ \\ 
\hline
\multirow{2}*{ED4} 	
&  Yellow & $2.00\pm 0.25$ & $0.99\pm 0.14$ & $2.61\pm 0.24$ & $2.17\pm 0.09$ \\
& Red & $1.45\pm 0.25$ & $0.75\pm 0.10$ & $0.51\pm 0.10$ & $2.05\pm 0.12$ \\ 
\hline
\multirow{2}*{ED5} 	
&  Yellow & $2.20\pm 0.26$ & $0.86\pm 0.14$ & $3.79\pm 09.33$ & $1.69\pm 0.16$ \\
& Red & $1.30\pm 0.20$ & $0.60\pm 0.11$ & $0.66\pm 0.13$ & $1.56\pm 0.14$ \\ 
\hline
\multirow{2}*{ED6} 	
&  Yellow & $1.08\pm 0.12$ & $1.14\pm 0.18$ & $2.67\pm 0.14$ & $1.76\pm 0.13$ \\
& Red & $0.53\pm 0.18$ & $0.70\pm 0.24$ & $0.67\pm 0.27$ & $1.39\pm 0.19$ \\ 
\hline
	\end{tabular}		
	\caption{Values of the NVA times $t^{\rm NVA/Y}_i$ and $t^{\rm NVA/R}_i$, $i=1, \ldots ,6$ corredsponding to the optimal point}
	\label{tab:app2}
\end{table}		

%
%\bibliographystyle{elsarticle-num}
%\bibliography{EDreference,simulation,derivativeFree,New_healthcare_reference,piermarini}

\end{document}